\documentstyle[amsfonts,12pt]{article}
\newcommand{\lon}{\longrightarrow}
\newcommand{\rar}{\rightarrow}
\newcommand{\hook}{\hookrightarrow}

\newcommand{\CP}{{\Bbb C} {\Bbb P}}

\newcommand{\p}{{\partial}}

\newcommand{\C}{{\Bbb C}}
\newcommand{\rH}{\mbox{H}}
\newcommand{\R}{{\Bbb R}}
\newcommand{\Z}{{\Bbb Z}}
\newcommand{\ot}{\otimes}

\newcommand{\Img}{\mbox{\rm Im}\,}
\renewcommand{\Re}{\mbox{\rm Re}\,}
\renewcommand{\div}{\mbox{\rm div}\,}

\newcommand{\End}{\mbox{\rm End}\,}
\newcommand{\Ber}{\mbox{\rm Ber}\,}
\newcommand{\Beq}{\begin{equation}}
\newcommand{\Eeq}{\end{equation}}
\newcommand{\Beqr}{\begin{eqnarray*}}
\newcommand{\Eeqr}{\end{eqnarray*}}
\newcommand{\Bi}{\begin{itemize}}
\newcommand{\Ei}{\end{itemize}}
\newcommand{\Ba}{\begin{array}}
\newcommand{\Ea}{\end{array}}
\newcommand{\f}{{\cal O}}
\newcommand{\cT}{{\cal T}}
\newcommand{\cM}{{\cal M}}
\newcommand{\cN}{{\cal N}}

\newcommand{\cU}{{\cal U}}
\newcommand{\cF}{{\cal F}}

\newcommand{\cV}{{\cal V}}
\newcommand{\cX}{{\cal X}}
\newcommand{\cZ}{{\cal Z}}
\newcommand{\cY}{{\cal Y}}

\newcommand{\al}{\alpha}
\newcommand{\be}{\beta}
\newcommand{\ga}{\gamma}
\newcommand{\var}{\varepsilon}

\newcommand{\om}{\omega}
\newcommand{\dal}{\dot{\alpha}}
\newcommand{\dbe}{\dot{\beta}}
\newcommand{\dga}{\dot{\gamma}}

\newcommand{\tOm}{\widehat{\Omega}}
\newcommand{\tom}{\widehat{\omega}}
\newcommand{\tcY}{\widehat{\cY}}
\newcommand{\wtheta}{\widehat{\theta}}

\newcommand{\sip}{\smallskip}
\newcommand{\bip}{\bigskip}

\begin{document}
\sloppy

 \title{The extended moduli space of special \\ Lagrangian submanifolds}

\author{S.A.\ Merkulov \\
{\small\em Department of Mathematics, University of Glasgow}}

\date{}
\maketitle
\sloppy

\begin{abstract}
It is well known that the moduli space of all deformations of a compact
special Lagrangian submanifold $X$ in a Calabi-Yau manifold $Y$ within
the class of special Lagrangian submanifolds is isomorphic to the first
de Rham cohomology group of $X$. Reinterpreting the embedding data $X\subset Y$
within the mathematical framework of the Batalin-Vilkovisky
quantization, we find a natural deformation problem which extends the above
moduli space to the full de Rham cohomology group of $X$.
\end{abstract}

\bip

\begin{center}
{\bf \S 1. Introduction}
\end{center}

Let $Y$ be a Calabi-Yau manifold of complex dimension $m$
with K\"ahler form $\om$ and a
nowhere vanishing holomorphic $m$-form $\Omega$. A
compact real $m$-dimensional submanifold $X\hook Y$ is called
{\em special Lagrangian}\, if
$\omega|_X=0$ and $\Img\,\Omega|_X=0$. According to McLean \cite{McL},
the moduli space of all deformations of $X$ inside $Y$
within the class of special Lagrangian submanifolds
is a smooth manifold whose
tangent space at $X$ is isomorphic to the first de Rham cohomology group $\rH^1(X,\R)$.

\sip

Moduli spaces of special Lagrangian submanifolds are playing an
increasingly important role in quantum cohomology and related topics. On
physical grounds, Strominger, Yau and Zaslow argued \cite{SYZ} that
whenever a Calabi-Yau manifold $Y$ has a mirror partner $\hat{Y}$, then
$Y$ admits a foliation $\rho: Y^{2m} \rar B^m$ by special Lagrangian
tori $T^m$ and $\hat{Y}$ is the compactification of the family of dual
tori $\hat{T}^m$ along the fibres of the projection $\rho$ (for the
mathematical account of this construction see \cite{Mor}).

\sip

Recently, the mirror conjecture has been extended by Vafa \cite{Vafa}, also on physical grounds,
to Calabi-Yau manifolds $\hat{Y}$ equipped with stable vector bundles $W$.
According to Vafa, the mirror partner of such a pair $(\hat{Y},W)$
must be a triple $(Y, X, L)$ consisting of a Calabi-Yau manifold $Y$, a
compact special Lagrangian submanifold $X\hook Y$ and a flat unitary
line bundle $L$ on $X$ together with an
isomorphism between the moduli space (with typical tangent space $\rH^1(\hat{Y}, \End
W)$) of all
deformations of the holomorphic vector bundle $W\rar \hat{Y}$, and the moduli space
(with typical tangent  space $\rH^1(X,\R)\ot \C$) associated with McLean's
deformations of the embedding $X\hook
Y$ and deformations of the flat unitary line bundle $L$ on $X$. Actually,
Vafa conjectures that much more must be true:
$$
\rH^*(\hat{Y}, \End W) = \rH^*(X,\R)\ot \C.
$$

\sip

This raises a problem of finding a geometric interpretation
of the full de Rham cohomology group  of a special Lagrangian submanifold
$X\hook Y$. \footnote{This also raises much easier problems of extending the
moduli space of holomorphic vector bundles $W$ on $\hat{Y}$ with typical tangent
space $\rH^1(\hat{Y},\End W)$ to a supermanifold with typical tangent superspace
$\rH^*(\hat{Y}, \End W)$, and extending the
moduli space of flat unitary line bundles on $X$ with typical
tangent space $\rH^1(X, \R)$
to a supermanifold with typical tangent superspace  $\rH^*(X, \R)$; the
 latter two extensions are very
straightforward within the Batalin-Vilkovisky formalism and will  be
discussed elsewhere.}
Its solution  is the main theme of the present
paper. By moving into the mathematical realm of Batalin-Vilkovisky
quantization, we devise, out of the {\em same}\, data $X\hook Y$,
a deformation problem whose moduli space has the typical tangent space
isomorphic to $\rH^*(X,\R)$ thereby extending McLean's moduli space to the
full de Rham group.  The idea is very simple.

\sip

First, out of $Y$ we
construct a real $(2m|2m)$-dimensional supermanifold $\cY:= \Pi \Omega^1 Y$,
where $\Omega^1 Y$ is the {\em real}\, cotangent bundle and $\Pi$
denotes the parity change functor. The
supermanifold $\cY$ comes equipped with a complex structure
and a nowhere vanishing holomorphic section, $\tOm\in\Gamma(\cY,
\Ber_{holo}(\cY))$, of the holomorphic Berezinian bundle
induced by the holomorphic $m$-form form $\Omega$ on $Y$. We
note that if $\cX\subset \cY$ is an $(m|m)$-dimensional  real slice,
then $\tOm$ restricts to $\cX$ as a global no-where vanishing section
of the bundle $\C\ot \Ber(\cX)$.

\sip

Second, we construct a real $(2m|2m+1)$-dimensional
supermanifold $\tcY:= \cY \times \R^{0|1}$ and note that it comes
canonically equipped with an odd exact contact structure represented by a
$1$-form $\wtheta$. We also note that the K\"ahler form $\om$ on $Y$ gives rise
to an even smooth function $\tom$ on $\tcY$ and call an $(m|m)$-dimensional
sub-supermanifold
$\cX\hook \tcY$ {\em special Legendrian}\, if the
following conditions hold
$$
\wtheta|_{\cX}=0, \ \ \tom|_{\cX}=0, \ \ \ \Img(\tOm|_{p(\cX)})=0,
$$
where $p$ denotes the natural projection $\tcY \rar \cY$.

\sip

If $\cX\hook \tcY$ is special Legendrian, then $\cX_{red}\hook Y$ is
special Lagrangian. If $X \hook Y$ is special Lagrangian (with normal bundle
denoted by $N$), then the
associated supermanifold $\cX:= \Pi N^*$ is a special Legendrian
sub-supermanifold of $\tcY$. However, the correspondence between the
special Lagrangian submanifolds $X\hook Y$ and special Legendrian
sub-supermanifolds $\cX\hook \tcY$ is {\em not}\, one-to-one ---
passing from the Calabi-Yau manifold $Y$ to the associated contact
supermanifold $\tcY$ brings precisely the right amount of new degrees of
freedom to extend McLean's moduli space $\rH^1(X,\R)$ to the full de Rham
group $\rH^*(X,\R)$.

\bip

\noindent{\bf Main theorem}. {\em Let $X \hook Y$ be a compact special
Lagrangian submanifold of a Calabi-Yau manifold and let $\cX=\Pi N^*\hook
\tcY$ be the associated special Legendrian sub-supermanifold of the contact
supermanifold. The maximal moduli space $\cM$ of  deformations of $\cX$
inside $\cY$ within the class of special Legendrian sub-supermanifolds is a
smooth supermanifold whose tangent superspace at $\cX$ is canonically isomorphic
to}\, $\Pi \rH^*(X, \R)$.

\bip

In view of Vafa's conjectures, it is important to study geometric
structures induced on $\cM$ from the original data $X\hook Y$ (cf.\
\cite{Hit,Tyu}). It should be also noted that

\sip

The paper is organised as follows. In \S2 and 3 we study
extended moduli spaces of general compact submanifolds $X$ in a
manifold $Y$. In \S 4 we specialize to the case when the ambient
manifold $Y$ has a symplectic structure and the submanifolds $X\hook Y$
are Lagrangian. In \S 5 we consider another special case when $Y$
is a complex manifold equipped with a nowhere-vanishing
holomorphic volume form $\Omega$ and $X\hook Y$ is a real slice
of $Y$ satisfying $\Img \Omega|_X=0$. Finally, in \S 6 we combine all
the previous results to prove the Main Theorem.

\bip

\bip

\begin{center}
{\bf \S 2. Extended Kodaira moduli spaces}
\end{center}

\bip

{\bf 2.1. Families of compact submanifolds.}
Let $Y$ and $M$ be smooth manifolds and let $\pi_{1}: Y\times  M\lon Y$
and $\pi_{2} : Y\times M \lon M$ be the natural projections. A{\em
family of compact submanifolds of the manifold $Y$ with the moduli
space $M$}\, is a  submanifold $F\hookrightarrow Y\times M$ such that
the restriction, $\nu$, of the projection $\pi_{2}$ on $F$ is a proper regular
map. Thus the family $F$ has the structure of a double fibration
$$
Y \stackrel{\mu}{\longleftarrow} F \stackrel{\nu}{\lon} M,
$$
where $\mu\equiv \pi_{1}\mid_{F}$. For each $t\in M$, there is an associated
compact submanifold $X_t$ in $Y$ which is said to belong to the
family $F$. Sometimes we use a more explicit notation $\{X_t\hook Y \mid t\in M\}$
to denote the family $F$ of compact submanifolds.
The family $F$ is called {\em maximal}\, if for any other  family
$\tilde{F}\hookrightarrow Y\times \tilde{M}$ such that
$\mu\circ\nu^{-1}(t) = \tilde{\mu}\circ\tilde{\nu}^{-1}(\tilde{t})$ for some
points
$t\in M$ and $\tilde{t}\in \tilde{M}$, there is a  neighbourhood
$\tilde{U}\subset \tilde{M}$ of the point $\tilde{t}$ and a smooth map
$f: \tilde{U} \rar M$ such that $\tilde{\mu}\circ\tilde{\nu}^{-1}(\tilde{t}') =
\mu\circ\nu^{-1}\left(f(\tilde{t}')\right)$ for every $\tilde{t}'\in \tilde{U}$.

\sip

Similar definitions can be made in the category of complex manifolds,
category of (complex) supermanifolds and category of analytic (super)spaces.

\bip

{\bf 2.2. The Kodaira map.} Consider a $1$-parameter family, $F\hook Y\times M$,
of compact sub(super)manifolds in a (super)manifold $Y$, where
$M=\R^{1|0}$ or $\R^{0|1}$ with the natural coordinate denoted by $t$ (such a
family is often called a {\em $1$-parameter deformation of the
sub(super)manifold $X=\mu\circ \nu^{-1}(0)$}).
There is a finite covering $\{U_i\}$ of $F$ such that
the restriction to each $U_i$ of the ideal sheaf $J_F$ of $F\hook Y\times M$
is finitely generated, say $J_F|_{U_i}=\langle f^{\al}_i \rangle$,
$\al=1, \ldots, \mbox{codim}\, F$. It is easy to see that the family
$\{ \frac{\p f^{\al}_i}{\p t}\bmod J_F\}$ defines a global section of the normal
bundle, $N_F$, of the embedding $F\hook Y\times M$ and hence gives rise
to a morphism of sheaves,
$$
\Ba{rccc}
k:& TM & \lon & \nu_*^{0} N_F\\
&\frac{\p}{\p t} & \lon & \{ \frac{\p f^{\al}_i}{\p t}\bmod J_F\}.
\Ea
$$
This morphism, or rather its restriction
$$
k_t: T_t M \lon \rH^0(X_t, N_t),
$$
where $N_t\simeq N_F|_{\nu^{-1}(t)}$ is the normal bundle of $X_t\hook
Y$,  is called the {\em Kodaira map}.

\sip

If $M$ is the 1-tuple point $\R[t]/t^2$, then the family $F$ is
called an {\em infinitesimal deformation}\, of $X=\mu\circ \nu^{-1}(0)$ in
$Y$. The Kodaira map establishes a one-to-one correspondence between all
possible infinitesimal deformations of $X$ inside $Y$ and the vector
superspace $\rH^0(X,N)$. Often in this paper we shall be interested in
deformations of $X$ inside $Y$ within a  class of {\em special}\, (say,
complex, Lagrangian, Legendrian, etc.) submanifolds. The associated
set of all possible infinitesimal deformations of $X$ is  a vector subspace of
$\rH^0(X,N)$ called the
{\em Zariski tangent space at $X$ to the moduli space of (special) compact
submanifolds}. Note that we do {\em not}\, require that any element  of the
Zariski tangent space
at $X$ necessary exponentiates to a genuine 1-parameter deformation of $X$.
Put another way, the Zariski tangent space to the moduli space $\cM$
makes sense even when $\cM$ does not exists as a smooth manifold!

\sip

{\bf 2.3.  From manifolds to supermanifolds.} Given a compact
submanifold, $X\hook Y$, of a smooth manifold $Y$. The associated exact
sequence
$$
0 \lon TX \lon TY|_X \lon N \lon 0
$$
implies the canonical map
$$
\Lambda^*TY|_X \lon \Lambda^* N \lon 0
$$
which in turn implies the canonical embedding, $\cX \hook \cY$, of the associated
supermanifolds $\cX:= (X, \Lambda^* N)$ and $\cY:=(Y, \Lambda^* TY)$
(cf. \cite{Sch1,Sch2}). In
more geometrical terms, $\cX\simeq \Pi N^*$, $\cY\simeq \Pi \Omega^1 Y$
and $\cX\hook \cY$ corresponds just to the natural inclusion, $\Pi N^*
\subset \Pi \Omega^1 M|_X$. The supermanifold $\cY$ comes equipped
canonically with
an even {\em Liouville}\, 1-form $\theta$ defined, in a natural local coordinate system
$(x^a, \psi_a\simeq \Pi\p/\p x^a)$ on $\cY$, as follows
$$
\theta = \sum_{a=1}^{n} dx^a \psi_a, \ \ \ \  \ n=\dim Y.
$$
The odd two-form
$$
\eta:= d\theta= -\sum_{a=1}^{n}dx^a\wedge d\psi_a
$$
is non-degenerate and hence equips $\cY$ with an odd symplectic
structure.

\sip

A $(p,n-p)$-dimensional sub-supermanifold $\cX\hook \cY$ is called {\em
Lagrangian}\, if $\eta_{\cX}=0$ (this implies, in particular, that
$\theta|_{\cX}$ is closed). It is called {\em exact Lagrangian}\, if
$\theta|_{\cX}$ is an exact $1$-form on $\cX$.

\bip

{\bf 2.3.1. Lemma.} {\em For any submanifold $X\hook Y$, the associated
sub-supermanifold $\cX\hook \cY$ is exact Lagrangian.}

\sip

{\em Proof}. Assume $\dim\, X=p$.
We can always choose a local coordinate system $(U, x^a)$ in a
tubular neighbourhood $U$ of (a part of) $X$ inside $Y$ in such a way that
$X\cap U=\{x^a=0, \  a=p+1, \ldots, n\}$. Then the normal bundle $N$ of
$X\hook Y$ is locally generated by $\p/\p x^a$ with $a=p+1, \ldots, n$.
Hence $\cX \hook \cY$ is locally given by the equations $x^a=0,
\psi_b=0$ where $a=p+1, \ldots,n$ and $b=1, \ldots, p$. It is now obvious
that $\theta|_{\cX}=0$. Finally, $\dim \cX=(p,n-p)$. $\Box$

\bip

{\bf 2.3.2. Remark.} It also follows from the above proof that, for any
submanifold $X\hook Y$ with the normal bundle $N_X$,
$\theta|_{\Pi N_X^*}=0$. It is not hard to check that the reverse is also true:
if $\cX\hook \cY$ is an $(p|n-p)$-dimensional sub-supermanifold
such that $\theta|_{\cX}=0$, then $\cX=\Pi N^*X$ for some submanifold
$X\hook Y$.

\bip

{\bf 2.4. The extended Kodaira map.}
Let $\cX$ be a Lagrangian sub-supermanifold of a supermanifold $\cY$
equipped with an odd symplectic structure $\eta$. Then, as usually, one
gets an odd isomorphism $j:\Omega^1 \cX \stackrel{\eta^{-1}}{\rar} \cN$,
where $\cN$ is the
normal bundle of $\cX\hook \cY$. In particular, there is a monomorphism
of sheaves,
$$
i: \f_{\cX}/\R \stackrel{j\circ d}{\lon} \cN,
$$
where $d$ is the exterior derivative.
\sip

Consider now a one (even or odd) parameter family
of compact exact Lagrangian sub-supermanifolds of the supermanifold
$\cY=\Pi \Omega^1 Y$,
i.e.\ a double fibration
$$
\cY \stackrel{\mu}{\longleftarrow} \cF \stackrel{\nu}{\lon} M,
$$
with $\nu$ being a proper submersion and
$\cX_t:=\mu\circ \nu^{-1}(t)$ being a compact exact Lagrangian sub-supermanifold of
$(\cY, \eta)$ for every $t\in M\subset \R^{1|0}$ or $\R^{0|1}$ .
\sip

{\bf 2.4.1. Lemma}.
{\em For the family $\{\cX_t \hook \cY | t\in M\}$ as above the Kodaira map
$k_t: T_tM \rar \mbox{\em H}^0(\cX_t, N_t)$ factors as follows}
$$
k_t: T_tM \stackrel{k'}{\lon} \rH^0(\cX_t, \f_{\cX_t})/\R
\stackrel{i}{\lon} \rH^0(\cX_t, \cN_t).
$$
\sip

{\em Proof.} Since $\mu^*(\eta)|_{X_t}=0$, we have
$$
\mu^*(\eta)= A \wedge dt
$$
for some 1-form $A$ on $\cF$ whose restriction to $\nu^{-1}(t)$
represents, under the isomorphism
$j:\Omega^1 \cX_t \stackrel{\eta^{-1}}{\rar} \cN_t$, the normal vector field
$k_t(\p/\p t)$. On the other hand,
$$
\mu^*(\theta)= \Psi dt + dB,
$$
for some smooth functions $\Psi$ and $B$ on $\cF$ with parities
$\tilde{\Psi}=\tilde{t}+1$ and $\tilde{B}=1$\footnote{Here and elsewhere $\tilde{}$\
stands for the parity of the kernel symbol.}. Thus
$A= d\Psi$ completing the proof. $\Box$

\bip

{\bf 2.4.2. Corollary.} {\em Let $\cX$ be a compact exact Lagrangian
sub-supermanifold of an odd symplectic supermanifold $\cY$. Then
the Zariski tangent space at $\cX$ to the moduli space of all
deformations of $\cX$ within the class of exact Lagrangian sub-supermanifolds
is isomorphic to $\mbox{\em H}^0(\cX, \f_{\cX})/\R$}.

\bip

\bip

{\bf 2.4.3. Definition.} If $\cY=\Pi \Omega^1 Y$ and
$\cX_t \simeq \Pi N_t^*$, where $N_t$ is the
normal bundle of some submanifold $X_t\hook Y$, then
$\rH^0(\cX_t, \f_{\cX})/\R \simeq \rH^0(X_t, \Lambda^* N_t)/\R$. The
associated map  $k': T_tM \lon \rH^0(X_t, \Lambda^* N_t)/\R$
is called the {\em extended Kodaira map}.

\bip

{\bf 2.5. Extended Kodaira moduli space.} Kodaira \cite{Ko}
proved that if $X\hook Y$ is a compact complex submanifold of a complex
manifold with $\rH^1(X,N)=0$, then there exists a maximal moduli space
$M$ parametrizing all possible deformations of $X$ inside $Y$ whose
tangent space at the point $X$ is isomorphic to $\rH^0(X,N)$.

\sip

With the same data $X\hook Y$ one associates a pair $\cX=\Pi N^* \hook
\cY=\Pi \Omega^1 Y$ and asks for all possible holomorphic deformations of $\cX$
inside $(\cY, \eta)$ within the class of complex exact Lagrangian sub-supermanifolds.

\bip

{\bf 2.5.1. Theorem.} {\em
Let $X \hook Y$ be a compact complex submanifold of a complex manifold
and $\cX\hook \cY$ the associated compact
complex exact Lagrangian sub-supermanifold. If
$H^1(X, \Lambda^k N)=0$ for all $k\geq 1$, then
there exists a maximal moduli space $\cM$, called the {\em extended
Kodaira moduli space}, which parametrizes all possible
deformations of $\cX$ inside $(\cY, \eta)$
within the class of complex exact Lagrangian sub-supermanifolds. Its tangent
space at the point $\cX$ is canonically isomorphic to $\sum_{k\geq 1} H^0(X,\Lambda^k N)$,
with the following $\Z_2$-grading:
$[\cT_{\cX}\cM]_0= \sum_{k\in 2\Z+1} H^0(X,\Lambda^k N)$
and $[\cT_{\cX}\cM]_1= \sum_{k\in 2\Z} H^0(X,\Lambda^k N)$.}

\bip

{\em Proof}\, is routine, cf.\ \cite{Ko,Me}.

\bip

{\bf 2.5.2. Example.} Let $X$ be a projective line $\CP^1$ embedded into
a complex 3-fold $Y$ with normal bundle $N=\f(1)\oplus \f(1)$. In this case
the Kodaira moduli space is a complex $4$-fold $M$ canonically equipped,
according to Penrose, with a self-dual conformal structure,
while the extended Kodaira moduli space $\cM$ is a $(4|3)$-dimensional
supermanifold isomorphic to $\Pi \Omega^2_+ M$, where $\Omega^2_+ M$ is
the bundle of self-dual 2-forms on $M$.

\bip

\bip

\begin{center}
{\bf \S 3. Restoring the lost constants}
\end{center}

\bip

{\bf 3.1. Odd contact structure.}
Let $X$ be a compact submanifold of a manifold $Y$. It is shown in \S 2
that the Zariski tangent space to the extended moduli space of
Lagrangian deformations of $\cX=\Pi N^*$ inside $\cY= \Pi \Omega^1 Y$ is
$\Pi H^0(X, \Lambda^* N)/\R$. One can easily restore the lost constants
$\R$ by extending $\cY$ to an odd contact supermanifold $\tcY$ and
studying {\em Legendrian}\,  families of compact sub-supermanifolds in
$\tcY$.

\sip

Consider $\tcY:= \cY\times \R^{0|1}$ and define a 1-form on $\tcY$
$$
\wtheta= d\var + p^*(\theta),
$$
where $p: \tcY \rar \cY$ is the natural projection and $\var$ is the
standard coordinate on $\R^{0|1}$.
The form $\wtheta$ defines an odd contact structure on $\tcY$.

\bip

{\bf 3.1.1. Lemma.} {\em For any submanifold $X\hook Y$, the associated
sub-supermanifold $\cX=\Pi N^* \hook \tcY$ is Legendrian with respect to the odd
contact structure $\theta$.}

\sip

{\em Proof.}  $\wtheta|_{\cX}= d\var|_{\cX}+ p^*(\theta)|_{\cX}=0+ 0=0$. $\Box$

\bip

Thus one can associate with data $X\hook Y$ the moduli space $\cM$ of all possible
deformations of $\cX$ inside $\tcY$ within the class of Legendrian
sub-supermanifolds.

\bip

{\bf 3.1.2. Proposition.} {\em The Zariski tangent space to $\cM$ at $\cX$ is}\,
$\Pi \rH^0(X, \Lambda^* N)$.

\sip

{\em Proof.} If
$
\tcY \stackrel{\hat{\mu}}{\longleftarrow} \cF \stackrel{\nu}{\lon}
M\subset \R^{1|0}\ \mbox{or}\ \R^{0|1}
$
is a 1-parameter family of compact Legendrian sub-supermanifolds, then
$$
\hat{\mu}^*(\wtheta)= \Psi dt
$$
for some $\Psi\in \Gamma(\cF, \f_{\cF})$ with $\tilde{\Psi}=\tilde{t}+1$
(compare this with
$\mu^*(\eta)=d\Psi\wedge dt$ in 2.4.1). The restriction of $\Psi$
to $\nu^{-1}(t)$ represents the image of $\p/\p t$ under the extended
Kodaira map. $\Box$

\bip

{\bf 3.2. Remark.} If $\{\cX_t\hook \tcY |\, t\in \cM\}$ is a family of
compact Legendrian sub-supermanifolds, then $\{p(\cX_t)\hook \cY |\, t\in \cM\}$
is a family of {\em exact}\,  Lagrangian sub-supermanifolds.
\sip


\bip

{\bf 3.3. An important observation.} If $(Y, \om)$ is a symplectic
manifold and $X\hook Y$ a Lagrangian submanifold with respect to $\om$,
then the normal bundle
$N$ is canonically isomorphic to $\Omega^1 X$. Thus  the associated
extended Zariski tangent space is isomorphic to $\Omega^*X$.

\bip

\begin{center}
{\bf \S 4. Even + odd symplectic structures}
\end{center}

{\bf 4.1. Isotropic Lagrangian sub-supermanifolds.}
In this section we assume that $Y$ is an even $2m$-dimensional symplectic
manifold. The symplectic 2-form $\omega$ on $Y$  gives rise to a global even
function $\tom$  on the associated odd symplectic supermanifold $\cY=
\Pi \Omega^1Y$ (and hence  on $\tcY= \cY \times \R^{0|1}$)
defined, in a natural local coordinate system $(x^a, \psi_a=\Pi \p/\p x^a)$,
as follows
$$
\tom= \sum_{a,b=1}^{2m} \omega^{ab}(x)\psi_a \psi_b,
$$
where $\omega^{ab}(x)$ is the matrix inverse to the matrix $\omega_{ab}(x)$ of
components of $\omega$ in the basis $dx^a$. The latter function gives
rise to an odd Hamiltonian vector field $Q$ on $\cY$ (or a contact
vector field $Q$ on $\hat{\cY}$) defined by
$$
Q \lrcorner\, \eta = d\tom,
$$
and is given, in a natural local coordinate system, by
$$
Q=\sum_{a,b}\omega^{ab}\psi_b \frac{\p}{\p x^a} + \sum_{a,b,c,d,e}
w^{ad}\frac{\p \om_{bc}}{\p x^a} \om^{ce}\psi_e\frac{\p}{\p \psi_b}.
$$

\sip

Differentials forms on $Y$ can be identified with smooth functions on the
supermanifold $\Pi TY$, $\Gamma(Y, \Omega^*Y)= \Gamma(\Pi TY, \f_{\Pi TY})$.
Under this identification the de Rham differential $d: \Omega^* Y \rar
\Omega^* Y$ corresponds to an odd vector field $d$ on $\Pi TY$
satisfying $d^2=0$. The even symplectic form $\om$ establishes an isomorphism
$\phi: \Pi TY \rar \cY$ and hence maps $d$ into an odd vector
field $\phi_*d$ on $\cY$ which, as it is not hard to check \cite{AKSZ},
coincides precisely with $Q$. This observation implies, in
particular, that $Q^2=0$ and $Q\tom=0$.

\bip

{\bf 4.1.1 Definition.} A Lagrangian (resp.\ Legendrian) sub-supermanifold
of $(\cY, \eta)$ (resp.\ $(\hat{\cY}, \wtheta)$) is called
$\om$-{\em isotropic}\, if $\tom|_{\cX}=0$ (resp.\ $\tom|_{p(\cX)}=0$).

\sip

If $\cX\hook \cY$ is $\om$-isotropic, then $Q|_{\cX}\in \Gamma(\cX,\cT\cX)$.

\bip

{\bf 4.1.2 Lemma.} {\em Let $(Y,\om)$ be a symplectic manifold and let
$\cY:= \Pi \Omega^1 Y$.

\sip

 (i) If $\cX$ is a compact  $(m|m)$-dimensional
$\om$-isotropic Lagrangian submanifold $(\cY, \eta)$, then $\cX_{red}$ is
a compact Lagrangian submanifold of $(Y,\om)$.

\sip

(ii) Let $X$ be a compact Lagrangian submanifold of $(Y, \om)$. Then the
associated compact $(m|m)$-dimensional sub-supermanifold
$\cX:=\Pi N^*_X \hook \cY$ is $\om$-isotropic. Moreover, under the isomorphism
$\Gamma(\cX, \f_{\cX})= \Gamma(X, \Omega^*X)$ the
vector field $Q|_{\cX}$ goes into the usual de Rham differential $d$ on $X$}.

\sip

{\em Proof}\, is very straightforward when one uses Darboux coordinates.

\bip

{\bf 4.2 Normal exponential map.} Let $X$ be an $r$-dimensional compact manifold
of an $n$-dimensional manifold $Y$ and let $\cX=\Pi N^*$ be the
associated  Lagrangian sub-supermanifold
of the odd symplectic supermanifold $(\cY=\Pi \Omega^1
Y, \eta)$.
\bip

{\bf 4.2.1. Lemma}. {\em  There exist
\Bi
\item a tubular neighbourhood $\cU$  of $\cX$ in $\cY$,
\item a tubular neighbourhood $\cV$ of $0_{\cX}$ in $\Pi \Omega^1
\cX$, where $0_{\cX}\simeq \cX$ is the zero section of the bundle
$\Pi \Omega^1 \cX \rar \cX$,
\item a diffeomorphism
$\exp: \cV \rar \cU$,
\Ei
such that
\Bi
\item[(i)] $\exp|_{0_{\cX}}: \cX\rar \cX$ is the identity map,  and
\item[(ii)] $\exp^*(\theta) - \theta_0 = d F$ for some $F\in \Gamma(\cX,
\f_{\cX})$,
\Ei
where $\theta$ is the Liouville form on $\Pi \Omega^1 Y$ and $\theta_0$ is
the Liouville form on $\Pi \Omega^1 \cX$. In particular,
$\exp^*(\eta)=\eta_0$, where $\eta_0$ is the natural odd symplectic
structure on $\Pi \Omega^1 \cX$. }
\bip

{\em Proof}.
There is a tubular neighbourhood $U$ of $\cX_{red}$ in $Y$ which can be
identified via the normal exponential map with a tubular neighbourhood
$V\subset N$ of
the zero section of the normal bundle $N$ of $X$ in $Y$. These
neighbourhoods and the exponential map have a canonical extension to the
map $\exp: {\cU}\rar {\cV}$ which has the property (i). We only have to
check the validity of (ii). Let $(x^{\al}, x^{\dal})$, $\al=1,\ldots,r$,
$\dal=r+1, \ldots,n$,  be a local
trivialisation of $N$, where $x^{\al}$ are local coordinates on the base
of $N$ and $x^{\dal}$,
are the fibre coordinates. In the associated local
coordinate system
$(x^{\al}, x^{\dal}, \psi_{\al}:= \Pi
\p/\p x^{\al}, \psi_{\dal}:=\Pi \p/\p x^{\dal})$ on $\cV\subset \Pi\Omega^1 \cX$  the zero
section $0_{\cX}$ is given by the equations $x^{\dal}= \psi_{\al}=0$.
We have
$$
\exp^*(\theta)= dx^{\al} \Pi\frac{\p}{\p x^{\al}} +
dx^{\dal} \Pi\frac{\p}{\p x^{\dal}}= dx^{\al} \psi_{\al} +
dx^{\dal}\psi_{\dal},
$$
and
$$
\theta_0  = dx^{\al} \Pi\frac{\p}{\p x^{\al}} +
d\psi_{\dal} \Pi\frac{\p}{\p \psi_{\dal}}= dx^{\al} \psi_{\al} -
d\psi_{\dal}x^{\dal}.
$$
Hence $\exp^*(\theta) - \theta_0 = d(\psi_{\dal}x^{\dal})$. Since
$\psi_{\dal}x^{\dal}$ is an invariant, the
statement follows. $\Box$

\bip

{\bf 4.2.2. Remark.}
The above Lemma establishes a one-to one correspondence between nearby to $\cX$ exact
Lagrangian sub-supermanifolds and global exact differential forms on
$\cX$. Note, however, that this correspondence is not canonical but
depends on the choice of the normal exponential map $\exp:\cV\rar \cU$.
If $f$ is a global odd section of $\f_{\cX}$ (such that $df\in \cV\subset
\Pi \Omega^1 \cX$) and $\cX_{df}\hook \cY$ is the
associated Lagrangian sub-supermanifold, then we have a diffeomorphism
$$
\exp_{df}: \cX \stackrel{df}{\lon} \cV \stackrel{\exp}{\lon} \cX_{df}.
$$

\sip

Consider now a particular case when $Y$ is a $2m$-dimensional symplectic manifold
$(Y, \om)$ and $X\hook Y$ is a compact Lagrangian submanifold with
respect to $\om$. In this case the normal bundle $N$ of $X$ in $Y$ is
isomorphic to $\Omega^1 X$ and hence the total space of $N$ is
naturally a symplectic manifold implying that $\Pi \Omega^1 \cX$ (with
$\cX=\Pi N^*$) comes canonically equipped with an odd vector field $Q_0$
such that $Q_0^2=0$ (cf. subsection 4.1). Since the normal exponential
map $N\supset V \stackrel{\exp}{\lon} U\subset Y$  can be
chosen to be a symplectomorphism, the associated extended exponential
map $\Pi \Omega^1 \cX \supset \cV \stackrel{\exp}{\lon} \cU\subset \Pi \Omega^1 Y$
can be chosen to satisfy the additional property
$$
\exp_*(Q_0)=Q.
$$
Note also that the isomorphism $N=\Omega^1 X$ implies $\cX=\Pi \cT X$ which
in turn implies $\Gamma(\cX, \f_{\cX})= \Omega^*X$. Then we have

\bip

{\bf 4.2.3. Lemma.} {\em For any $(\cU, \cV, \exp)$ as above and any
exact Lagrangian submanifold $\cX_{df}\hook \cU$,
the function $\exp_{df}^*(\tom|_{\cX_{df}})\in   \Gamma(\cX,
\f_{\cX})= \Omega^*X$ defines a closed (non-homogeneous, in
general) differential form on $X$.}

\sip

{\em Proof}. Let $\phi_{df}: \Pi \Omega^1 \cX \rar \Pi \Omega^1 \cX$ be
a translation by $df$ along the fibres of the projection $\Pi \Omega^1
\cX\rar \cX$. In the natural coordinates on $\Pi\Omega^1 \cX$ we have
\Beqr
\left[Q_0 - (\phi_{df})_*Q_0 \right]\exp^*(\tom) & = &
\left[\sum_{\al}\left(\frac{\p f}{\p x^{\al}} + \sum_{\dbe}\psi_{\dbe}
\frac{\p^2 f}{\p x^{\dbe} \p \psi_{\dal}}\right)\frac{\p}{\p
x^{\dal}}\right.\\
&& - \left. \sum_{\al,\dbe} \psi_{\dbe} \frac{\p^2 f}{\p x^{\al} \p
x^{\be}}\frac{\p}{\p \psi_{\al}}
\right] \sum_{\ga}\psi_{\dga}\psi_{\ga}\\
&=& \sum_{\dal,\dbe} \psi_{\dal} \psi_{\dbe} \frac{\p^2 f}{\p x^{\al}
\p x^{\be}}\\
&=& 0,
\Eeqr
and hence
\Beqr
Q_0|_{\cX}\left(\exp_{df}^*(\tom|_{\cX_{df}})\right)&=
&Q_0|_{\cX}\left(\phi_{df}^*\circ\exp^*(\tom)\right)|_{\cX}\\
&=& \left[\phi_{df}^*\circ\left(Q_0\exp^*(\tom)\right)\right]|_{\cX}\\
&=& \left[\phi_{df}^*\circ\exp^*(Q\tom)\right]|_{\cX}\\
&=&0.
\Eeqr
Then the statement follows from the fact that under the isomorphism $
\Gamma(\cX,\f_{\cX})= \Omega^*X$ the vector field $Q_0|_{\cX}\in
\Gamma(\cX, \cT\cX)$ goes into
the de Rham differential on $\Omega^* X$ (cf.\ Lemma 4.1.2(ii)). $\Box$

\bip

{\bf 4.3. Moduli space of isotropic sub-supermanifolds.}
Given a compact Lagrangian submanifold $X$
of a symplectic manifold $(Y,\om)$.
With these data one may associate
the extended moduli space $\cM$ of all possible deformations of $\cX=\Pi N^*$
inside $\hat{\cY}$ within the class of Legendrian, $\om$-isotropic
sub-supermanifolds.

\sip

{\bf 4.3.1. Theorem.} {\em The Zariski tangent space to $\cM$ is
$\Pi \mbox{\em H}^0(X, \Omega^*X_{closed})$, where $\Omega^*X_{closed}$ is the
sheaf
of closed differential forms on $X$}.

\sip

{\em Proof}. Let $\{\cX_t\hook \tcY |t\in M\}$ be a
1-parameter family of $\om$-isotropic Legendrian submanifolds, and let
$$
\cY \stackrel{\mu}{\longleftarrow} \cF \stackrel{\nu}{\lon} M,
$$
be an associated 1-parameter family of Lagrangian, $\om$-isotropic
sub-supermanifolds. The vector field $V_f$ on $\cY$ gives rise
to a vector field on $\cY\times M$ (denoted by the same letter)
which is tangent to $\cF \hook \cY \times M$. We have
$$
V_f|_{\cF}\, \lrcorner\, \mu^*(\eta) = (V_f|_{\cF} \,\lrcorner\, d\Psi)\wedge dt,
$$
implying
$$
\mu^*(df) = (V_f|_{\cF}\Psi)dt.
$$
Since $\mu^*(df)=d(f|_{\cF})=0$, we get $V_f|_{\cF}\Psi=0$. Finally, the
required statement follows 4.1.1(ii) which says that $V_f|_{\cF}$ is
essentially the de Rham differential. $\Box$

\bip

\bip

\begin{center}
{\bf \S 5. Moduli spaces of special real slices}
\end{center}
{\bf 5.1. Batalin-Vilkovisky structures.} Let $\cY$ be an $(n|n)$-dimensional
compact supermanifold equipped with an odd symplectic form $\eta$ and an even
nowhere-vanishing section $\mu$ of the Berezinian bundle $\Ber(\cY)$.
Such data have been extensively studied by
A.S.\ Schwarz in \cite{Sch1,Sch2} in the context of Batalin-Vilkovisky
quantization.

\sip

The volume form $\mu$ induces the Berezin integral,
$\int_{\mu} f$,
 on smooth functions $f$ on $\cY$. In particular, $\mu$ gives rise to a divergence
operator $\div  V$ on smooth vector fields $V$ on $\cY$ which can
be characterized by the formula \cite{Ge}
$$
\int (\div V)f \mu =- \int V(f) \mu.
$$
If $x^a$ is a local coordinate system on $\cY$ and  $D^*(dx^a)$ the
associated local basis  of $\Ber(\cY)$, then $\mu=\rho(x) D^*(dx^a)$
for some even nowhere-vanishing even function $\rho(x)$ and
$$
\div V = \frac{1}{\rho}(-1)^{\tilde{a}(1+\tilde{V})}\frac{\p (V^a
\rho)}{\p x^a},
$$
where $\tilde{a}$ is the parity of $x^a$ and $V^a$ are the components of
$V$ in the basis $\p/\p x^a$. Another possible definition,
which also works in the holomorphic category, is
$$
\div V= \frac{L_V \mu}{\mu},
$$
where $L_V$ stands for the Lie derivative along the vector field $V$.

\sip

In particular, if $V_f$ is the hamiltonian vector field on $\cY$
associated to a smooth function $f\in \Gamma(\cY, \f_{\cY})$, then one
defines a second order operator,
$$
\Delta f := \frac{1}{2} \div V_f.
$$
Note that this operator depends solely on $\mu$ and $\eta$. The
situations when $\Delta^2=0$ are of special interest in the context of
Batalin-Vilkovisky quantization. The data
$(\cY, \eta, \mu)$ with property $\Delta^2=0$ are sometimes called
$SP$-{\em manifolds} \cite{Sch1,Sch2} or {\em Batalin-Vilkovisky supermanifolds}
\cite{Ge}.
The structure $(\cY, \eta, \mu)$ which arises in the context
of Calabi-Yau manifolds does actually satisfy the requirement
$\Delta^2=0$, see \S 6.

\bip

{\bf 5.2. Integration on Lagrangian sub-supermanifolds.}
Let $\cY$ again be an $(n|n)$-dimensional
compact oriented supermanifold equipped with an odd symplectic form $\eta$ and an even
nowhere-vanishing section $\Omega$ of the Berezinian bundle $\Ber(\cY)$,
and let $\cX\hook \cY$ be a compact $(r|n-r)$-dimensional Lagrangian
sub-supermanifold. Then the extension
$$
0\lon \cT \cX \lon \cT \cY|_{\cX} \lon \cN \lon 0
$$
and the isomorphism $\cN\simeq \Pi \Omega^1 \cX$ imply
$$
\Ber (\cY)|_{\cX} = \Ber(\cX)^{\ot 2}.
$$
Thus the volume form $\tOm$ on $\cY$ induces a volume form on $\cX$ which
we denote by $\tOm^{1/2}$. A possible problem with taking the
square root is overcome with the assumption that $\cY$ is
oriented;  a clear and explicit construction of
$\tOm^{1/2}$ is given by A.S.\ Schwarz in \cite{Sch1}.

\sip

As an example, let us consider the case when $\cY=\Pi \Omega^1 Y$, where
$Y$ is an $n$-dimensional compact manifold equipped with a
nowhere-vanishing $n$-form $\Omega$. The latter gives rise, via the isomorphism
$\Ber(\cY)\simeq\mbox{Det}\, (Y)^{\ot 2}$, to a volume form $\tOm$ on
$\cY$. If $x^{a}$ is a local coordinates system on $Y$ in which $\Omega =
\al(x) dx^1\wedge \ldots \wedge dx^n$, then, in the associated
local coordinate system $(x^a, \psi_a:=\Pi \p/\p x^a)$ on $\cY$,
$$
\tOm=\al^2(x) D^*(dx^a, d\psi_a).
$$
In particular, if $\cX\hook \cY$ is a Lagrangian sub-supermanifold given
locally by the equations $x^b=0, \psi_e=0$, $b=r+1,\ldots,n$,
$e=1,\ldots,r$, then $\tOm^{1/2}= \al(x)|_{\cX_{red}}D^*(dx^e,
d\psi_b)$.

\sip

There is a natural morphism of sheaves,
$$
\Ba{rccc}
F: & \Omega^*Y & \lon & \f_{\cY}\\
&\sum
w_{a_1\ldots a_k} dx^{a_1}\wedge\ldots \wedge dx^{a_k} & \lon &
\sum \al(x)^{-1} w_{a_1\ldots a_k} \var^{a_1 \ldots a_k a_{k+1}\ldots
a_n}\psi_{a_{k+1}}\ldots \psi_{a_n},
\Ea
$$
where $\var^{a_1 \ldots a_k a_{n}}$ is the antisymmetric tensor with
$\var^{1\ldots n}=1$. One has \cite{Wi,Sch1},
$$
F(dw)=\Delta F(w),
$$
where $\Delta$ is the Batalin-Vilkovisky operator on $(\cY, \tOm,
\eta)$.

\bip

{\bf 5.2.1. Lemma} {\em Let $Y$ be a
manifold $Y$ equipped with a nowhere vanishing volume form $\Omega$ and
let $\cY=\Pi \Omega^1 Y$. If, for any compact submanifold $X\hook Y$,
the function  $\Phi\in \Gamma(\cY, \f_{\cY})$ is such that
$$
\int_{\Pi N^*_X}\Phi\,{\tOm^{1/2}}  =0,
$$
then $\Phi=\Delta \Psi$ for some
$\Phi\in \Gamma(\cY, \f_{\cY})$}.

\bip

{\em Proof}. We may assume for simplicity that $\Phi$ is homogeneous in odd
coordinates $\psi_a$, i.e. that $\Phi=F(w)$ for some $k$-form on $Y$.
According to A.S.\ Schwarz \cite{Sch1},
$$
\int_{\Pi N^*_X} \Phi \,\tOm^{1/2}  = \int_X w.
$$
Since this vanishes for any compact submanifold $X\hook Y$, $w=ds$ for
some $(k-1)$-form $s$ on $Y$. Then $\Phi=F(w)=F(ds)=\Delta F(s)$. $\Box$

\bip

{\bf 5.3. Holomorphic volume forms.}
Let $Y$ be an $m$-dimensional complex manifold equipped with a no-where vanishing
holomorphic $m$-form $\Omega$. Then the associated $(m|m)$-dimensional complex
supermanifold\footnote{In this and the next
sections the subscript $c$ is used to distinguish holomorphic objects
from the real ones. In particular, $\Omega^1_c Y$ denotes the bundle of
holomorphic 1-forms on $Y$ as opposite to $\Omega^1 Y$ which denotes the
bundle of real smooth 1-forms on the real manifold underlying $Y$.} $\cY= \Pi
\Omega^1_{c}Y$ comes equipped with two
natural odd symplectic structures. The first one is holomorphic and is represented,
in a natural local coordinate system $(z^{\al},
\zeta_{\al}:= i\Pi\p/\p z^{\al})$, $\al=1, \ldots,m$, by the odd
holomorphic 2-form
$$
\eta_{c} = d\left(\sum_{\al} dz^{\al}\zeta_{\al}\right)=-\sum_{\al}d(x^{\al} + i
x^{\dal})\wedge d(\psi_{\dal} + i
\psi_{\al}),
$$
where $z^{\al}=x^{\al} + i x^{\dal}$, $\psi_{\al}=\Pi \p/\p x^{\al}$
and $\psi_{\dal}=\Pi \p/\p x^{\\dal}$. The second one is real and comes
from the identification of the real $(2m|2m)$-dimensional supermanifold
underlying $\cY$ (which we denote by the same letter $\cY$)
 with the real cotangent bundle $\Pi\Omega^1 Y$. It is
given by
$$
\eta= d\sum_{\al}\left( dx^{\al}d\psi_{\al} +
dx^{\dal}\psi_{\dal}\right).
$$
Clearly, $\eta=\Img \eta_{c}$.

\sip

The holomorphic $m$-form $\Omega$ induces, via the isomorphism
$$
\Ber_{c}(\cY)= [\Omega^m_{c} Y]^{\ot 2},
$$
a holomorphic volume form $\tOm$ on $\cY$.

\sip

A compact real $(m|m)$-dimensional sub-supermanifold $\cX \hook \cY$
is called a {\em real slice}\, if the sheaf $\C\ot\cT \cX$ is
isomorphic to the sheaf of smooth sections of $\cT_c \cY|_{\cX}$.
In this case $\tOm$ induces \cite{AKSZ} a smooth section, $\tOm|_{\cX}$, of the
complexified Berezinian bundle $\C\ot \Ber(\cX)$. A real slice $\cX\hook
\cY$ is called {\em special}\, if $\Img (\tOm|_{\cX})=0$. In this case
$\Re (\tOm|_{\cX})$ is a nowhere-vanishing real volume form on $\cX$.
If $\cX\hook \cY$ is also Lagrangian with respect to the real odd
symplectic structure $\eta$, then $\eta_c|_{\cX}$ is non-degenerate and hence
makes $\cX$ into an
odd symplectic manifold. According to 5.1, the data $(\Re (\tOm|_{\cX},
\eta_c|_{\cX})$ induces on the structure sheaf of $\cX$
a second-order differential operator $\Delta$.
Note that if $X\hook Y$ is real slice of the manifold $Y$ such that
$\Img \Omega|_X=0$, then the associated sub-supermanifold $\cX:=\Pi N^*
\hook \cY$ is a special Lagrangian real slice.

\bip

{\bf 5.3.1. Theorem} {\em Let $Y$ be a
complex manifold $Y$ equipped with a nowhere-vanishing holomorphic
$m$-form $\Omega$, $X$  a compact real slice of $Y$ such that
$\Img\Omega|_X=0$, and $\cX=\Pi N^*$ the associated special Lagrangian real
slice in $\cY=\Pi \Omega^1 Y$. Then the Zariski tangent space to the
moduli space of all possible deformations of $\cX$ inside $\cY$ within
the class of special Lagrangian real slices is isomorphic to the kernel
of the operator $\Delta: \Gamma(\cX, \f_{\cX})/\R \rar \Gamma(\cX,
\f_{\cX})$.
}

\bip

{\em Proof}. Consider a 1-parameter family, $\{ \cX_t\hook \cY\, |\,
t\in \R^{1|0}\ \mbox{or}\ \R^{0|1}\}$, of special Lagrangian real slices
in $\cY$ such that $\cX_{t=0}=\cX$. Let $z^{\al}=x^{\al}+ i x^{\dal}$ be a local coordinate
system on $Y$ in which $X$ is given by $x^{\dal}=0$. Then in the
associated local coordinate system $(z^{\al},
\zeta_{\al}:= i\Pi\p/\p z^{\al}=\psi_{\dal}+i\psi_{\al})$ on $\cY$, the equations
of $\cX_t\hook \cY$ are
$$
x^{\dal}= \frac{\p \Phi}{\p \psi_{\dal}}, \ \ \ \psi_{\al}= - \frac{\p
\Phi}{\p x^{a}}
$$
for some 1-parameter family of smooth functions $\Phi=\Phi(x^{\al},
\psi_{\dal}, t)$ satisfying the boundary condition $\Phi(x^{\al},
\psi_{\dal}, 0)=0$. The image of $\p/\p t$ under the extended Kodaira
map is represented by the function (see Lemma 2.4.1)
$$
k_{t=0}\left(\frac{\p}{\p t}\right)= \left.\frac{\p \Phi}{\p
t}\right|_{t=0}
\equiv \Psi.
$$
\sip

If $\tOm= \rho(z^{\al}, \zeta_{\al})D^*(dz^{\al}, d\zeta_{\al})$
for some holomorphic function $\rho(z^{\al}, \zeta_{\al})$, then
\Beqr
\tOm|_{\cX_t} & = & \rho\left(x^{\al}+ i\frac{\p \Phi}{\p \psi_{\dal}},
\psi_{\dal} - i\frac{\p \Phi}{\p x_{\al}}\right) D^*
\left(d(x^{\al}+ i\frac{\p \Phi}{\p \psi_{\dal}}),
d(\psi_{\dal} - i\frac{\p \Phi}{\p x_{\al}})\right) \\
&=&\rho\left(x^{\al}+ i\frac{\p \Phi}{\p \psi_{\dal}},
\psi_{\dal} - i\frac{\p \Phi}{\p x_{\al}}\right)
\Ber \left(\Ba{cc} \delta^{\al}_{\be} + i \frac{\p^2 \Phi}{\p x^{\be}\p
\psi_{\dal}} & \frac{\p^2 \Phi}{\p x^{\al}\p x^{\be}} \\
\frac{\p^2 \Phi}{\p \psi^{\dbe}\p \psi^{\dal}} &
\delta^{\al}_{\be} - i \frac{\p^2 \Phi}{\p x^{\al}\p
\psi_{\dbe}} \Ea \right) D^*(dx^{\al}, d\psi_{\dal}).
\Eeqr
Hence,
\Beqr
\left.\frac{d \Img(\tOm|_{\cX_t})}{dt}\right|_{t=0} &=&
\frac{1}{\rho_0}\sum_{\al}\left(\frac{\p \Psi}{\p \psi_{\dal}}
\frac{\p \rho_0}{\p x^{\al}} - \frac{\p \Psi}{\p x^{\al}}
\frac{\p \rho_0}{\p \psi_{\dal}} + 2\rho_0 \frac{\p^2 \Psi}{\p x^{\al}\p
\psi_{\dal}}\right)\rho_0 D^*(dx^{\al}, d\psi_{\dal})\\
&=& \left(\div H_{\Psi}\right) \Re \tOm|_{\cX}\\
&=& \left(\Delta \Psi\right) \Re \tOm|_{\cX},
\Eeqr
where $\rho_0=\Re\rho(x_{\al},\psi_{\dal})$. Hence $\Delta \Psi=0$.
$\Box$

\bip

\bip

\begin{center}
{\bf \S 6. Existence of the extended moduli space \\
 of special Lagrangian submanifolds}
\end{center}

{\bf 6.1. Initial data.}
Let $X$ be a compact special
Lagrangian submanifold of a  Calabi-Yau manifold $Y$ equipped with the K\"ahler
form $\om$ and a holomorphic volume form $\Omega$, and let $\cX=\Pi N^*\hook
\tcY$ be the associated special Legendrian sub-supermanifold of the contact
supermanifold $\tcY$ (see \S 1). With these data one naturally
associates the moduli superspace $\cM$ of all deformations of $\cX$
inside $\tcY$ within the class of special Legendrian sub-supermanifolds.

\bip

{\bf 6.2. Proposition}. {\em The Zariski tangent superspace to $\cM$ at $\cX$
is canonically isomorphic to $\Pi \mbox{\em H}^*(X, \R)$.}

\sip

{\em Proof}. It is not hard to check that under the isomorphism
$\f_{\cX}=\Omega^* X$ the Batalin-Vilkovisky operator $\Delta:
\f_{\cX}\rar \f_{\cX}$ goes into $2*d*$, where $d$ is the de Rham
differential and $*$ is the Hodge duality operator. Then Theorems 4.3.1 and
 5.3.1 imply that the Zariski tangent superspace is isomorphic to
$$
\Pi \Gamma(X, \Omega^*X_{closed}) \cap \Pi \Gamma(X, \Omega^*X_{coclosed})
= \Pi \rH^*(X,\R). \ \ \ \ \ \ \ \ \  \Box
$$

\bip

{\bf 6.3. Theorem. } {\em $\cM$ is a smooth supermanifold}.

\bip

{\em Proof}\, (after McLean \cite{McL}).
 Let $\cV$ be a tubular neighbourhood of the zero section
in $\Pi \Omega^1 \cX$, $\cU$ a tubular neighbourhood of $\cX$ in $\tcY$
and $\exp:\cV \rar \cU$ the normal exponential map constructed as in section 4.2.
This map identifies nearby (to $\cX$) special Legendrian
sub-supermanifolds $\cX_f$ of $\tcY$ with global odd sections $f$
of $\Gamma(\cX,\f_{\cX})$ and induces a diffeomorphism $\exp_f: \cX \rar
\cX_f$.
Let $\cV'$ be an open subset in $\Gamma(\cX,\f_{\cX})$ lying in the
preimage of $\cV$ under the map $d: \f_{\cX}\rar \Omega^1\cX$. We define
a non-linear map
$$
\phi: \cV'\subset \Gamma(\cX,\f_{\cX}) \lon \Omega^*X \bigoplus \Omega^*X
$$
as follows
$$
\phi(f)= \left(\exp_{f}^*(\tom),\ \
\Img\left(\frac{(p\circ\exp_{f})^*(
\tOm|_{p(\cX_f})}{\Re (\tOm|_{p(\cX)})}\right)^{1/2}\right),
$$
where $p: \tcY \rar \cY=\Pi \Omega^1 Y$ is the natural projection. The
square root in the above formula always exists (cf.\ section 5.2).
Note that $\phi^{-1}(0,0)=\cM$.

\sip

It follows from Lemma~4.2.3 that $\exp_{f}^*(\tom)\in \Omega^*X$ is a
closed differential form. Replacing $f$ with $tf$, we see that the map
$\exp_f:\cX \rar \tcY$ is homotopic to the inclusion $\cX\rar
\tcY$. Therefore, denoting by $[\ ]$ the cohomology class, we get
$[\exp_{f}^*(\tom)]=[\tom|_{\cX}]=0$ and conclude that
$\exp_{f}^*(\tom)$ is an exact differential form on $X$.

\sip

Since $\tOm$ is holomorphic, the
integral $\int_{p(\cX_f)}\tOm|_{p(\cX_f)}$
depends only on the homology class of $\cX_{red}$ in $Y$ \cite{AKSZ}.
Analogously, for any compact $(r|m-r)$-dimensional Lagrangian\footnote{with
respect to the odd symplectic structure induced on $p(\cX_f)$ from the
holomorphic odd symplectic structure on $\cY$, see section 5.2}
sub-supermanifold $\cZ\subset \cX_f$, the integral $\int_{\cZ}
\tOm^{1/2}$ (and hence its real and imaginary parts) depends only on
the homology class of $\cZ_{red}$ in $Y$. Since $\cZ_{red}$ is
homologous to an $r$-dimensional cycle in $\cX$ and $\Img
(\tOm|_{\cX}^{1/2})$ vanishes, we conclude that
$\int_{\cZ} \Img (\tOm^{1/2})=0$ for any such $\cZ$.
Thus, for any  smooth cycle $Z\hook X \subset Y$,
we have
\Beqr
0&=& \int_{\cZ_f} \Img (\tOm^{1/2})\\
&=& \int_{\Pi N^*_Z} \left(\frac{(p\circ\exp_{f})^*
\Img(\tOm|_{\cZ_f}^{1/2})}
{\Re (\tOm|_{\Pi N^*_Z})^{1/2}}\right) \Re (\tOm|_{\Pi N^*_Z})^{1/2}\\
&=& \int_{\Pi N^*_Z}
\Img\left(\frac{(p\circ\exp_{f})^*(
\tOm|_{\cZ_f})}{\Re (\tOm|_{\Pi N^*_Z})}\right)^{1/2} \Re
(\tOm|_{\Pi N^*_Z})^{1/2},
\Eeqr
where $\cZ_f:=p\circ \exp_f(\Pi N^*_Z)$ and we used the fact that
 $Z$ and $(\cZ_f)_{red}$ are homologous in $Y$.
By Lemma~5.2.1 and the fact that in our case $\Delta=*d*$,
the integrand of the last integral
is a coexact differential form in $\Omega^*X$.

\sip

Thus we proved that $\phi$ maps $\cV'\subset \Omega^*X$ into the subset
$$
\Omega^*X_{exact} \bigoplus \Omega^* X_{coexact} \subset
\Omega^*X\bigoplus \Omega^*X.
$$
Put another way, as a map from $C^{1,\al}$ differential forms on $X$
to exact and coexact $C^{0,\al}$ differential forms, $\phi$ is
surjective. Then, by
the Banach space implicit function theorem and elliptic regularity,
the extended moduli space
$\cM=\phi^{-1}(0,0)$ is smooth with tangent space at $0$ canonically isomorphic to
the kernel of the following operator (see the proofs of Theorems~4.3.1 and
5.3.1),
$$
\left.\frac{d}{dt} \phi(tf)\right|_{t=0} = (d, *d*) :\ \Omega^*X \lon
\Omega^*X \bigoplus \Omega^*X,
$$
which is precisely $\Pi \rH^*(X,\R)$. $\Box$

\bip

{\em Acknowledgement}. It is a pleasure to thank A.N.\ Tyurin for
valuable discussions.

\pagebreak

\bip

\noindent\mbox{\small Department of Mathematics, Glasgow University}

\noindent\mbox{\small 15 University Gardens, Glasgow G12 8QW, UK}

\noindent\mbox{\small e-mail: sm@maths.gla.ac.uk}


\begin{thebibliography}{99}

\bibitem{AKSZ} M.\ Alexandrov,\ M.\ Kontsevich, A.\ Schwarz and O.\
Zabolonsky, {\em The geometry of the master equation and topological
quantum field theory}, Int.\ J.\ Mod.\ Phys. {\bf A12} (1997), 1405-1430; hep-th/9502010.

\bibitem{Ge} E.\ Getzler, {\em Batalin-Vilkovisky algebras and
two-dimensional topological field theories}, Commun.\ Math.\ Phys.
{\bf 159} (1994), 265-285;  hep-th/9212043.


\bibitem{Hit} N.J.\ Hitchin,  {\em The moduli space of special
Lagrangian submanifolds}, dg-ga/9711002.

\bibitem{Ko} K.\ Kodaira, {\em A theorem of completeness of
characteristic systems for analytic families of compact submanifolds of
complex manifolds},  Ann.\ Math. {\bf  75} (1962), 146-162.

\bibitem{McL} R.C.\ McLean, {\em Deformations of calibrated submanifolds},
Duke University preprint, January 1996.

\bibitem{Me} S.A.\ Merkulov, {\em Existence and geometry of Legendre
moduli spaces}, Math.\ Z. \ {\bf 226 } (1997), 211-265.

\bibitem{Mor} D.R.\ Morrison, {\em The geometry underlying mirror
symmetry}, In Proc.\ European Algebraic Geometry Conf. (Warwick, 1996);
alg-geom/9608006.

\bibitem{Sch1} A.\ Schwarz, {\em Geometry of Batalin-Vilkovisky
quantization}, Commun.\ Math.\ Phys. {\bf 155} (1993) 249-260; hep-th 9205088.

\bibitem{Sch2} A.\ Schwarz, {\em Semiclassical approximation in
Batalin-Vilkovisky formalism}, Commun.\ Math.\ Phys. {\bf 158} (1994) 265-285;
hep-th/9210115.

\bibitem{SYZ} A.\ Strominger, S.-T.\ Yau and E.\ Zaslow, {\em Mirror
symmetry is $T$-duality}, Nucl.\ Phys. B {\bf 479} (1996), 243-259.

\bibitem{Tyu} A.N.\ Tyurin, {\em Special Lagrangian geometry and
slightly deformed algebraic geometry}, math.AG/9806006.

\bibitem{Wi} E.\ Witten, {\em A note on the antibracket formalism}, Mod.\
Phys.\ Lett.\ {\bf A5} (1990), 487.

\bibitem{Vafa} C.\ Vafa, {\em Extending mirror conjecture to Calabi-Yau
with bundles}, hep-th/9804131.


\end{thebibliography}
\end{document}